\documentclass[12pt,reqno]{amsart}
\usepackage{amsmath, amsthm, amscd, amsfonts, amssymb, graphicx, color}
\usepackage[bookmarksnumbered, colorlinks, plainpages]{hyperref}
\theoremstyle{definition}

\theoremstyle{remark}

\numberwithin{equation}{section}
\begin{document}
\begin{center}
{\textbf{Conformal Einstein soliton within the framework of para-K\"{a}hler manifold }}
\end{center}
\vskip 0.3cm
\begin{center}By\end{center}\vskip 0.3cm
\begin{center}
{Soumendu Roy \footnote{The first author is the corresponding author, supported by Swami Vivekananda Merit Cum Means Scholarship, Government of West Bengal, India.}, Santu Dey $^2$ and Arindam~~Bhattacharyya $^3$}
\end{center}
\vskip 0.3cm
\address[Soumendu Roy]{Department of Mathematics,Jadavpur University, Kolkata-700032, India}
\email{soumendu1103mtma@gmail.com}

\address[Santu Dey]{Department of Mathematics, Bidhan Chandra College, Asansol - 4, West Bengal-713303 , India}
\email{santu.mathju@gmail.com}

\address[Arindam Bhattacharyya]{Department of Mathematics,Jadavpur University, Kolkata-700032, India}
\email{bhattachar1968@yahoo.co.in}
\vskip 0.5cm
\begin{center}
\textbf{Abstract}\end{center}
The object of the present paper is to study some properties of para-K\"{a}hler manifold whose metric is conformal Einstein soliton. We have studied some certain curvature properties of para-K\"{a}hler manifold admitting conformal Einstein soliton.
\vspace {0.3cm}
\section{\textbf{Introduction}}
The notion of Einstein soliton was introduced by G. Catino and L. Mazzieri \cite{CATINO} in 2016, which generates self-similar solutions to Einstein flow,
\begin{equation}\label{1.1}
  \frac{\partial g}{\partial t} = -2(S-\frac{r}{2}g),
\end{equation}
where $S$ is Ricci tensor, $g$ is Riemannian metric and $r$ is the scalar curvature.\\
The equation of the Einstein soliton \cite{blaga} is given by,
\begin{equation}\label{1.2}
  \pounds_V g+2S+(2\lambda-r)g=0,
\end{equation}
where $\pounds_V$ is the Lie derivative along the vector field $V$, $S$ is the Ricci tensor, $r$ is the scalar curvature of the Riemannian metric $g$, and $\lambda$ is a real constant.\\
In 2015, N. Basu and A. Bhattacharyya \cite{nbab} introduced the notion of Conformal Ricci soliton equation \cite{soumendu}, \cite{roy}, given by \\
\begin{equation}\label{1.3}
\pounds_V g + 2S = [2\lambda - (p + \frac{2}{n})]g,
\end{equation}\\
 where $\pounds_V$ is the Lie derivative along the vector field $V$, $S$ is the Ricci tensor, $\lambda$ is constant, $p$ is a scalar non-dynamical field(time dependent scalar field) and $n$ is the dimension of the manifold. \\
So we introduce the notion of Conformal Einstein soliton as:
\begin{equation}\label{1.4}
   \pounds_V g+2S+[2\lambda-r+(p+\frac{2}{n})]g=0,
\end{equation}
where  $\pounds_V$ is the Lie derivative along the vector field $V$ , $S$ is the Ricci tensor, $r$ is the scalar curvature of the Riemannian metric $g$, $\lambda$ is real constant, $p$ is a scalar non-dynamical field(time dependent scalar field)and $n$ is the dimension of manifold.\\
Also it is said to be shrinking, steady or expanding according as $\lambda < 0$, $\lambda = 0$ and $\lambda > 0$ respectively.\\\\
In the present paper we study conformal Einstein soliton on para-K\"{a}hler manifold. The paper is organized as follows:\\
After introduction, section 2 is devoted for preliminaries on $n$- dimensional para-K\"{a}hler manifold, where $n$ is even. In section 3, we have studied conformal Einstein soliton on para-K\"{a}hler manifold. Here we proved if a $n$- dimensional para-K\"{a}hler manifold admits conformal Einstein soliton then the vector field associated with the soliton is solenoidal depens on the scalar curvature. We have also characterized the nature of the manifold if the manifold is quasi conformally flat,  pseudo-projectively flat and $W_2$- flat.\\
\vspace {0.3cm}
\section{\textbf{Preliminaries}}
Let $M$ be a connected differentiable manifold of dimension $n = 2m$, $m\geq2$, $F$ be a (1,1)-tensor field and $g$ be a pseudo-Riemannian metric on $M$. Then $(M, F, g)$ is said to be a para-K\"{a}hler manifold if the following conditions hold:
\begin{equation}\label{2.1}
  F^2=I, \quad g(FX,FY)=-g(X,Y),\quad \nabla F=0.
\end{equation}
for any $X,Y \in \chi(M)$, being the Lie algebra of vector fields on $M$, $\nabla$ is the Levi-Civita connection of $g$ and $I$ is the identity operator.\\
In a para-K\"{a}hler manifold $(M, F, g)$, the Riemannian curvature tensor $R$, the Ricci tensor $S$ and the scalar curvature $r$ are defined by:
\begin{equation}
  \widetilde{R}(X,Y,Z,W)=g(R(X,Y)Z,W), \nonumber
\end{equation}
\begin{equation}
  S(X,Y)= trace\{Z\rightarrow R(Z,X)Y\}, \nonumber
\end{equation}
\begin{equation}
  r=trace S.\nonumber
\end{equation}
Also the following properties are satisfied in a para-K\"{a}hler manifold:
\begin{equation}\label{2.2}
  R(FX,FY)Z=-R(X,Y)Z,
\end{equation}
\begin{equation}\label{2.3}
  R(FX,Y)Z=-R(X,FY)Z,
\end{equation}
\begin{equation}\label{2.4}
  S(FX,Y)=-S(FY,X),
\end{equation}
\begin{equation}\label{2.5}
 S(FX,FY)=-S(X,Y).
\end{equation}
\vspace {0.3cm}
\section{\textbf{Some results on conformal Einstein soliton within the framework of para-K\"{a}hler manifold }}
From equation \eqref{1.4}, we can write,
\begin{equation}\label{3.1}
  (\pounds_V g)(X,Y)+2S(X,Y)+[2\lambda-r+(p+\frac{2}{n})]g(X,Y)=0,
\end{equation}
for any $X,Y \in \chi(M)$, being the Lie algebra of vector fields on $M$.\\
Taking $X = e_i$, $Y = e_i$ in the above equation and summing over $i=1,2,....,n$, we get,
\begin{equation}\label{3.2}
  div V+r+[\lambda-\frac{r}{2}+\frac{1}{2}(p+\frac{2}{n})]n=0.
\end{equation}
If $V$ is solenoidal then $div V = 0$ and so from the above equation we have $r = \frac{2\lambda n}{n-2}+\frac{n(p+\frac{2}{n})}{n-2}$. Again if $r = \frac{2\lambda n}{n-2}+\frac{n(p+\frac{2}{n})}{n-2}$ then \eqref{3.2} gives $div V =0.$\\
So we can state the following theorem:\\\\
\textbf{Theorem 3.1.} {\em If the metric of an $n$-dimensional para-K\"{a}hler manifold satisfies a conformal Einstein soliton then the vector field associated with the soliton is solenoidal iff the scalar curvature is $\frac{2\lambda n}{n-2}+\frac{n(p+\frac{2}{n})}{n-2}$, provided $n>2$.}\\\\
The notion of quasi-conformal curvature tensor was introduced by Yano and Sawaki \cite{yano} and it is defined by:
\begin{eqnarray}\label{3.3}
  C(X,Y)Z &=& \alpha R(X,Y)Z+ \beta [S(Y,Z)X-S(X,Z)Y+g(Y,Z)QX \nonumber \\
          &-& g(X,Z)QY]-\frac{r}{n}(\frac{\alpha}{n-1}+2\beta)[g(Y,Z)X-g(X,Z)Y],
\end{eqnarray}
where $\alpha$, $\beta$  are constants, $Q$ is the Ricci operator, defined by $g(QX , Y) = S(X,Y)$ and $n$ is the dimension of the manifold.\\
Moreover if $\alpha = 1$ and $\beta = -\frac{1}{n-2}$,the above equation reduces to conformal curvature tensor \cite{eisenhart}.\\
Again a manifold $(M^n , g)$ where $n>3$, is said to be quasi conformally flat if $C=0.$\\
Now in an $n$- dimensional para-K\"{a}hler manifold, we can define the Ricci tensor $S$ as:
\begin{equation}\label{3.4}
S(X,Y)=\frac{1}{2}\sum_{i=1}^{n}\epsilon_i  \widetilde{R}(e_i,F e_i,X,FY),
\end{equation}
where $\{e_1 , e_2 , .... , e_n\}$ is an orthonormal frame and $\epsilon_i$ is the indicator of $e_i$, $\epsilon_i = g(e_i , e_i)=1.$\\
Taking inner product in \eqref{3.3} by $W$, we get,
\begin{eqnarray}\label{3.5}
  g(C(X,Y)Z,W) &=& \alpha \widetilde{R}(X,Y,Z,W)+ \beta [S(Y,Z)g(X,W)-S(X,Z)g(Y,W) \nonumber \\
          &+& g(Y,Z)S(X,W)-g(X,Z)S(Y,W)] \nonumber\\
          &-&\frac{r}{n}(\frac{\alpha}{n-1}+2\beta)[g(Y,Z)g(X,W)-g(X,Z)g(Y,W)].
\end{eqnarray}
Now if the manifold is quasi-conformally flat then the above equation reduces to,
\begin{multline}\label{3.6}
 \alpha \widetilde{R}(X,Y,Z,W)+ \beta [S(Y,Z)g(X,W)-S(X,Z)g(Y,W) \\
  +g(Y,Z)S(X,W)-g(X,Z)S(Y,W)]-\frac{r}{n}(\frac{\alpha}{n-1}+2\beta)[g(Y,Z)g(X,W)\\
  -g(X,Z)g(Y,W)]=0.
\end{multline}
Putting $X = e_i, Y= Fe_i, W= FW$ in the above equation and summing over $i=1,2,....,n$ and also using \eqref{3.4}, \eqref{2.4}, we get,
\begin{multline}\label{3.7}
  2\alpha S(Z,W)+2\beta [S(Z,W)-S(FZ,FW)] \\
  -\frac{r}{n}(\frac{\alpha}{n-1}+2\beta)[g(Z,W)-g(FZ,FW)]=0.
\end{multline}
Again using \eqref{2.5} and \eqref{2.1} in the above equation, we obtain,
\begin{equation}\label{3.8}
   2\alpha S(Z,W)+4\beta S(Z,W)-\frac{2r}{n}(\frac{\alpha}{n-1}+2\beta)g(Z,W)=0,
\end{equation}
which reduces to,
\begin{equation}\label{3.9}
  (\alpha+2\beta)S(Z,W)=\frac{r}{n}(\frac{\alpha}{n-1}+2\beta)g(Z,W).
\end{equation}
Taking $Z= e_i, W= e_i$ in the above equation and summing over $i= 1,2,.....,n$, we have,
\begin{equation}\label{3.10}
   r=0,
\end{equation}
since $\alpha \neq 0$.\\
Now if $r = 0$, then from \eqref{3.9}, we get $S = 0$, provided $\alpha + 2\beta \neq 0$ i.e the manifold is locally flat if $\alpha + 2\beta \neq 0$.\\
Then \eqref{1.4} becomes,
\begin{equation}\label{3.11}
  (\pounds_V g)(X,Y)+[2\lambda+(p+\frac{2}{n})]g(X,Y)=0
\end{equation}
for any $X,Y \in \chi(M)$, being the Lie algebra of vector fields on $M$.\\
Putting $X =Y= e_i$ and summing over $i=1,2,.....,n$, we get,
\begin{equation}\label{3.12}
  \sum_{i=1}^{n}\epsilon_i (\pounds_V g)(e_i,e_i)+ \sum_{i=1}^{n}\epsilon_i [2\lambda+(p+\frac{2}{n})]g(e_i,e_i)=0,
\end{equation}
which reduces to,
\begin{equation}\label{3.13}
  div V+[\lambda+\frac{1}{2}(p+\frac{2}{n})]n=0.
\end{equation}
If $V$ is solenoidal then $div V = 0$ and so from \eqref{3.13} we have $\lambda+\frac{1}{2}(p+\frac{2}{n}) = 0$. Also if $\lambda+\frac{1}{2}(p+\frac{2}{n}) = 0$, \eqref{3.13} reduces to $div V = 0$.\\
This leads to the following:\\\\
\textbf{Theorem 3.2.} {\em If the metric of an $n$-dimensional quasi-conformally flat para-K\"{a}hler manifold satisfies a conformal Einstein soliton then the vector field associated with the soliton is solenoidal iff $\lambda+\frac{1}{2}(p+\frac{2}{n}) = 0$.}\\\\
The pseudo-projective curvature tensor $\overline{P}$ \cite{prasad}  is given by:
\begin{eqnarray}\label{3.14}
 \overline{P}(X,Y))Z &=& aR(X,Y)Z+b[S(Y,Z)X-S(X,Z)Y] \nonumber \\
                                &-& \frac{r}{n}(\frac{a}{n-1}+b)[g(Y,Z)X-g(X,Z)Y],
\end{eqnarray}
where $a, b \neq 0$ are constants.\\
Also a manifold $(M^n , g)$, is said to be pseudo-projectively flat if $\overline{P}=0.$\\
Taking inner product in \eqref{3.14} by $W$, we get,
\begin{eqnarray}\label{3.15}
 g(\overline{P}(X,Y))Z,W) &=& a\widetilde{R}(X,Y,Z,W)+b[S(Y,Z)g(X,W)-S(X,Z)g(Y,W)] \nonumber \\
                                &-& \frac{r}{n}(\frac{a}{n-1}+b)[g(Y,Z)g(X,W)-g(X,Z)g(Y,W)].\nonumber\\
\end{eqnarray}
Now if the manifold is pseudo-projectively flat then the above equation reduces to,
\begin{multline}\label{3.16}
 a\widetilde{R}(X,Y,Z,W)+b[S(Y,Z)g(X,W)-S(X,Z)g(Y,W)] \\
  -\frac{r}{n}(\frac{a}{n-1}+b)[g(Y,Z)g(X,W)-g(X,Z)g(Y,W)]=0.
\end{multline}
Putting $X = e_i, Y= Fe_i, W= FW$ in the above equation and summing over $i=1,2,....,n$ and also using \eqref{3.4}, \eqref{2.4}, we get,
\begin{equation}\label{3.17}
  2aS(Z,W)+b[S(Z,W)-S(FZ,FW)]-\frac{r}{n}(\frac{a}{n-1}+b)[g(Z,W)-g(FZ,FW)=0.
\end{equation}
Again using \eqref{2.5} and \eqref{2.1} in the above equation, we obtain,
\begin{equation}\label{3.18}
  (a+b)S(Z,W)-\frac{r}{n}(\frac{a}{n-1}+b)g(Z,W)=0.
\end{equation}
Taking $Z= e_i, W= e_i$ in the above equation and summing over $i= 1,2,..,n$, we have,
\begin{equation}\label{3.19}
   r=0,
\end{equation}
since $a \neq 0$.\\
Now if $r = 0$, then from \eqref{3.18}, we get $S = 0$, provided $a+b \neq 0$ i.e the manifold is locally flat if $a+b \neq 0$.\\
Then \eqref{1.4} becomes,
\begin{equation}\label{3.20}
  (\pounds_V g)(X,Y)+[2\lambda+(p+\frac{2}{n})]g(X,Y)=0
\end{equation}
for any $X,Y \in \chi(M)$, being the Lie algebra of vector fields on $M$.\\
Putting $X =Y= e_i$ and summing over $i=1,2,.,n$, we get,
\begin{equation}\label{3.21}
  \sum_{i=1}^{n}\epsilon_i (\pounds_V g)(e_i,e_i)+ \sum_{i=1}^{n}\epsilon_i [2\lambda+(p+\frac{2}{n})]g(e_i,e_i)=0,
\end{equation}
which reduces to,
\begin{equation}\label{3.22}
  div V+[\lambda+\frac{1}{2}(p+\frac{2}{n})]n=0.
\end{equation}
If $V$ is solenoidal then $div V = 0$ and so from \eqref{3.22} we have $\lambda+\frac{1}{2}(p+\frac{2}{n}) = 0$. Also if $\lambda+\frac{1}{2}(p+\frac{2}{n}) = 0$, \eqref{3.22} reduces to $div V = 0$.\\
This leads to the following:\\\\
\textbf{Theorem 3.3.} {\em If the metric of an $n$-dimensional pseudo-projectively flat para-K\"{a}hler manifold satisfies a conformal Einstein soliton then the vector field associated with the soliton is solenoidal iff $\lambda+\frac{1}{2}(p+\frac{2}{n}) = 0$.}\\\\
The $W_2$-curvature tensor $(n > 2)$ \cite{mishra} is given by:
\begin{equation}\label{3.23}
   W_2(X,Y)Z=R(X,Y)Z+\frac{1}{n-1}[g(X,Z)QY-g(Y,Z)QX].
 \end{equation}
Moreover a manifold is $W_2$- flat if $\widetilde{W_2}(X,Y,Z,U)=g(W_2(X,Y)Z,U)=0$.\\
Taking inner product in \eqref{3.23} by $U$, we get,
\begin{equation}\label{3.24}
   g(W_2(X,Y)Z,U)=\widetilde{R}(X,Y,Z,U)+\frac{1}{n-1}[g(X,Z)S(Y,U)-g(Y,Z)S(X,U)].
 \end{equation}
 Now if the manifold is $W_2$- flat then the above equation reduces to,
 \begin{equation}\label{3.25}
   \widetilde{R}(X,Y,Z,U)+\frac{1}{n-1}[g(X,Z)S(Y,U)-g(Y,Z)S(X,U)]=0.
 \end{equation}
 Putting $X = e_i, Y= Fe_i, U= FU$ in the above equation and summing over $i=1,2,....,n$ and also using \eqref{3.4}, \eqref{2.4}, we get,
 \begin{equation}\label{3.26}
   S(Z,U)+\frac{1}{n-1}[S(FZ,FU)-S(Z,U)]=0.
 \end{equation}
 Again using \eqref{2.5} in the above equation, we obtain,
 \begin{equation}\label{3.27}
   (n-3)S(Z,U)=0.
 \end{equation}
 Then we have $S(Z,U) = 0$,  for any $Z,U \in \chi(M)$, being the Lie algebra of vector fields on $M$, since $n$ is even.\\
 Hence from the above equation we can get, $r = 0.$\\
 Then \eqref{1.4} becomes,
\begin{equation}\label{3.28}
  (\pounds_V g)(X,Y)+[2\lambda+(p+\frac{2}{n})]g(X,Y)=0
\end{equation}
for any $X,Y \in \chi(M)$, being the Lie algebra of vector fields on $M$.\\
Putting $X =Y= e_i$ and summing over $i=1,2,.....,n$, we get,
\begin{equation}\label{3.29}
  \sum_{i=1}^{n}\epsilon_i (\pounds_V g)(e_i,e_i)+ \sum_{i=1}^{n}\epsilon_i [2\lambda+(p+\frac{2}{n})]g(e_i,e_i)=0,
\end{equation}
which reduces to,
\begin{equation}\label{3.30}
  div V+[\lambda+\frac{1}{2}(p+\frac{2}{n})]n=0.
\end{equation}
If $V$ is solenoidal then $div V = 0$ and so from \eqref{3.30} we have $\lambda+\frac{1}{2}(p+\frac{2}{n}) = 0$. Also if $\lambda+\frac{1}{2}(p+\frac{2}{n}) = 0$, \eqref{3.30} reduces to $div V = 0$.\\
This leads to the following:\\\\
\textbf{Theorem 3.4.} {\em If the metric of an $n$-dimensional $W_2$- flat para-K\"{a}hler manifold satisfies a conformal Einstein soliton then the vector field associated with the soliton is solenoidal iff $\lambda+\frac{1}{2}(p+\frac{2}{n}) = 0$.}\\\\


\begin{thebibliography}{}
\bibitem{nbab}Nirabhra Basu, Arindam Bhattacharyya, {\em Conformal Ricci soliton in Kenmotsu manifold}, Global
Journal of Advanced Research on Classical and Modern Geometries(2015), Vol. 4, Isu. 1, pp. 15-21.
\bibitem{blaga} A. M. Blaga, {\em On Gradient $\eta$-Einstein Solitons}, Kragujevac Journal of Mathematics(2018), Volume 42(2), pp-229–237.
\bibitem{CATINO} G. Catino and L. Mazzieri, {\em Gradient Einstein solitons}, Nonlinear Anal(2016). 132, pp-66–94.
\bibitem{eisenhart} L. P. Eisenhart, {\em Riemannian Geometry}, Princeton University Press, Princeton, N,J., (1949).
\bibitem{mishra} G. P. Pokhariyal, R. S. Mishra, {\em The curvature tensor and their relativistic significance},
Yokohoma Math. J. (1970), vol-18, pp- 105–108.
\bibitem{prasad}  B. Prasad, {\em A pseudo projective curvature tensor on a Riemannian manifold}, Bull. Calcutta
Math. Soc.(2002), 94 , pp- 163–166.
\bibitem{soumendu} Soumendu Roy, Arindam~~Bhattacharyya, {\em Conformal Ricci solitons on 3-dimensional trans-Sasakian manifold}, Jordan Journal of Mathematics and Statistics (2020), Vol- 13(1), pp-89-109.
\bibitem{roy} Soumendu Roy, Santu Dey, Arindam Bhattacharyya, Shyamal Kumar Hui, {\em $*$-Conformal $\eta$-Ricci Soliton on Sasakian manifold}, 	arXiv:1909.01318v1 [math.DG].
\bibitem{yano} K. Yano and S. Sawaki, {\em Riemannian manifolds admitting a conformal transformation group},
J. Differential Geometry(1968), 2 , pp-161-184.
\end{thebibliography}
\end{document}